\documentclass[12pt]{article}
\usepackage[cp1251]{inputenc}
\usepackage{amsmath}
\usepackage{amssymb}
\usepackage[T2A]{fontenc}
\usepackage{epsfig}
\usepackage{graphicx}

\begin{document}

\begin{large}

\centerline{\bf On the optimality of one integral inequality for
closed curves in $\mathbb R^4$ }

\bigskip

\centerline{ V. Gorkavyy, \,R. Posylaieva }

\end{large}
\bigskip

\bigskip

\begin{small}
{\it Abstract}. The optimality of the integral inequality
$\int\limits_\gamma\sqrt{k_1^2+k_2^2+k_3^2}ds>2\pi$ for closed
curves with non-vanishing curvatures in $\mathbb R^4$ is
discussed. We prove that an arbitrary closed curve of constant
positive curvatures in $\mathbb R^4$ satisfies the inequality
$\int\limits_\gamma\sqrt{k_1^2+k_2^2+k_3^2}ds \geq 2\sqrt{5}\pi$.


\medskip


\end{small}

\bigskip

\bigskip

{\bf \large Introduction}

\medskip

The famous Fenchel-Borsuk theorem of the classical theory of
curves states that the total curvature of an arbitrary smooth
closed curve $\gamma$ in $\mathbb R^n$ is greater than or equal to
$2\pi$:
\begin{equation}\label{Fenchel}
\int\limits_\gamma k_1\, ds\,\, \geq\,\, 2\pi,
\end{equation}
and the equality holds if and only if $\gamma$ is a convex closed
curve in $\mathbb R^2$, see [1,Ch.21], [2]-[3].

In [4] the first author described a series of integral
inequalities for curvatures of smooth closed curves in $\mathbb
R^n$ which may be viewed as a direct generalization of the
Fenchel-Borsuk inequality. Namely, let $\gamma$ be an arbitrary
smooth closed curve in $\mathbb R^n$, $n\geq 4$. Suppose that
$\gamma$ has non-vanishing curvatures $k_1$, $k_2$, ..., $k_j$ for
some $2\leq j\leq n-1$. Then the following inequality holds:
\begin{equation}\label{FenchelN}
\int\limits_\gamma \sqrt{k_{j-1}^2+k_{j}^2+k_{j+1}^2} ds > 2\pi.
\end{equation}
Consequently, if all the curvatures $k_1$, ..., $k_{n-1}$ of
$\gamma\subset \mathbb R^n$ are non-vanishing, then
(\ref{FenchelN}) holds true for each $2\leq j\leq n-1$, and thus
$\gamma$ satisfies a sequence of $n-2$ different integral
inequalities.

The inequality (\ref{FenchelN}) is {\it optimal} in the case of
odd $j$, see [4]. Actually, for any fixed odd $2\leq j\leq n-1$
one can construct a sequence of smooth closed curves
$\left\{\gamma_m \right\}_{m=1}^\infty$ in $\mathbb R^n$ so that
the values of $\int\limits_{\gamma_m}
\sqrt{k_{j-1}^2+k_{j}^2+k_{j+1}^2} ds$ tend to $2\pi$ as
$m\to\infty$. If $n$ is even, then the desired sequence
$\left\{\gamma_m \right\}_{m=1}^\infty$ may consists of closed
curves of {\it constant} curvatures in $\mathbb R^n$; if $n$ is
odd, then $\left\{\gamma_m \right\}_{m=1}^\infty$ in $\mathbb R^n$
may be obtained by slight perturbations of curves of constant
curvatures in $\mathbb R^{n-1}\subset \mathbb R^n$. Thus curves of
constant curvatures provide the optimality for (\ref{FenchelN}) in
the case of odd $j$.

As for the case of even $j$, the problem of the optimality of
(\ref{FenchelN}) still remains a quite challenging and interesting
open problem, which motivated this research paper.

We would start to discuss the problem by considering the simplest
case of $n=4$. As stated above, an arbitrary smooth closed curve
$\gamma\subset \mathbb R^4$ satisfies two inequalities:
\begin{equation}\label{4_1}
\int\limits_\gamma \sqrt{k_{1}^2+k_{2}^2+k_{3}^2} ds > 2\pi,\quad
\hbox{if }\, k_1>0, k_2>0;
\end{equation}
\begin{equation}\label{4_2}
\int\limits_\gamma \sqrt{k_{2}^2+k_{3}^2} ds > 2\pi,\quad \hbox{if
}\, k_1>0, k_2>0, k_3>0.
\end{equation}

The inequality (\ref{4_2}) is optimal since it corresponds to the
odd value $j=3$, see [4].

As for the inequality (\ref{4_1}), it looks rather trivial in view
of (\ref{Fenchel}), and hence the intuition tells us that
(\ref{4_1}) is not optimal. The main result of the paper confirms
partially this supposition.

\medskip

{\bf Theorem 1.} {\it Let $\gamma$ be a smooth closed curve in
$\mathbb R^4$ with non-vanishing constant curvatures $k_1$, $k_2$
and $k_3$. Then the following optimal estimate holds:}
\begin{equation}\label{4_3}
\int\limits_\gamma \sqrt{k_{1}^2+k_{2}^2+k_{3}^2} ds \geq
2\sqrt{5}\pi.
\end{equation}

A computer added numerical analysis demonstrates that (\ref{4_3})
holds true for some closed curves with non-constant curvatures
too. Together with Theorem 1, this suggest us conjecturing that
the inequality (\ref{4_3}) remains true for any closed curve with
non-vanishing curvatures in $\mathbb R^4$.

Let us recall the genesis of the inequality (\ref{4_1}). For an
arbitrary smooth curve $\gamma\subset \mathbb R^4$ with
non-vanishing curvatures $k_1$ and $k_2$ one can consider a
well-defined family of two-dimensional osculating planes of
$\gamma$, which are spanned by the first and second vectors of the
Frenet frame of $\gamma$. This family of planes may be interpreted
as a smooth closed curve in the Grassmann manifold $G_{2,4}$, it
is called the osculating indicatrix of $\gamma$ and denoted by
$\tilde \gamma$. The Grassmann manifold $G_{2,4}$ is embedded into
$\mathbb R^6$ with the help of Plucker coordinates, see [5,
Ch.8.2], [6], and therefore $\tilde \gamma$ may be viewed as a
smooth closed curve in $\mathbb R^6$. It turns out that the total
first curvature of $\tilde \gamma$ is equal exactly to the
left-hand side of (\ref{4_1}), and thus the inequality (\ref{4_1})
for $\gamma$ is just the Fenchel-Borsuk inequality for
$\tilde\gamma$, c.f. [4].

If $\gamma$ has constant curvatures, then the the following
stronger result holds.

\medskip

{\bf Theorem 2.} {\it Let $\gamma$ be a smooth closed curve in
$\mathbb R^4$ with non-vanishing constant curvatures $k_1$, $k_2$
and $k_3$. Let $\tilde \gamma\subset \mathbb R^6$ be the
osculating indicatrix of $\gamma$. The one has the following:

1) $\tilde \gamma$ is a smooth closed curve;

2) $\tilde \gamma$ belongs to a four-dimensional affine plane\,
$\mathbb R^4\subset \mathbb R^6$;

3) the curvatures $\tilde k_1$, $\tilde k_2$ and $\tilde k_3$ of
$\tilde \gamma$ are non-vanishing constant;}

4) \begin{equation}\label{4_4} \int\limits_{\tilde \gamma} \tilde
k_{1} d\tilde s = \int\limits_\gamma
\sqrt{k_{1}^2+k_{2}^2+k_{3}^2} ds;
\end{equation}

5) \begin{equation}\label{equality} \int\limits_{\tilde \gamma}
\sqrt{\tilde k_{1}^2+ \tilde k_{2}^2+ \tilde k_{3}^2} d\tilde s =
\sqrt{2} \int\limits_\gamma \sqrt{k_{1}^2+k_{2}^2+k_{3}^2} ds.
\end{equation}

\medskip

It was very surprising for us to discover the relationship
(\ref{equality}) which looks quite elegant although individual
expressions of $\tilde k_j$ in terms of $k_i$ are rather
cumbersome.

Due to Theorem 2, the described procedure of constructing the
osculating indicatrix may be viewed as a particular transformation
of closed curves of constant curvatures in $\mathbb R^4$. The
transformation may be iterated, and then at every step we obtain a
new curve of constant curvatures in $\mathbb R^4$. This results in
a specific sequence of curves of constant curvatures in $\mathbb
R^4$ which is generated by the initial curve of constant curvature
$\gamma$. Remark that the value of $\int\limits_\gamma
\sqrt{k_{1}^2+k_{2}^2+k_{3}^2} ds$ is augmented by $\sqrt{2}$ at
each step of iteration.

It would be interesting to extend the proved theorems to
particular and/or general families of closed curves with {\it
non-constant} curvatures in $\mathbb R^4$.

\bigskip

\bigskip
\newpage

{\bf \large 1. Closed curves of constant curvatures in $\mathbb
R^4$}

\medskip

 Let $\gamma$ be a smooth curve in $\mathbb R^4$,
whose curvatures $k_1$, $k_2$, $k_3$ are constant and
non-vanishing. Then $\gamma$ may be parameterized as follows:
\begin{eqnarray}\label{position_vector} x_1 = a_1 \cos \alpha_1 t,\, x_2 = a_1 \sin \alpha_1 t,\,
 x_3 = a_2 \cos \alpha_2 t,\, x_4 = a_2 \cos \alpha_2 t,\end{eqnarray}
where $a_1$, $a_2$, $\alpha_1$ and $\alpha_2$ are constants, see
[1, Ch.33].

For an arc-length $s$ of $\gamma$ one has
\begin{eqnarray}\label{length}
\frac{ds}{dt} \, = \,
\sqrt{a_1^2\alpha_1^2+a_2^2\alpha_2^2},\end{eqnarray} hence one
needs to assume $a_1^2\alpha_1^2+a_2^2\alpha_2^2\not = 0$ to
guarantee the smoothness (regularity) of $\gamma$.

For calculating the curvatures of $\gamma$, one may apply standard
formulae of the classical theory of curves, see [1, Ch.32].
Elementary differential-geometric calculations result in the
following statement.

{\bf Proposition 3.} {\it The curvatures of $\gamma$ are expressed as follows:
\begin{eqnarray}\label{k1}
k_1 = \frac{\sqrt{a_1^2\alpha_1^4+a_2^2\alpha_2^4}}{a_1^2\alpha_1^2+a_2^2\alpha_2^2},
\end{eqnarray}
\begin{eqnarray}\label{k2}
k_2 = \frac{a_1a_2\alpha_1\alpha_2\vert\alpha_1^2-\alpha_2^2\vert}
{(a_1^2\alpha_1^2+a_2^2\alpha_2^2)\sqrt{a_1^2\alpha_1^4+a_2^2\alpha_2^4}},
\end{eqnarray}
\begin{eqnarray}\label{k3} k_3 = \frac{\alpha_1\alpha_2}{\sqrt{a_1^2\alpha_1^4+a_2^2\alpha_2^4}}.
\end{eqnarray}}

Therefore, in order to guarantee the smoothness of $\gamma$ and
the non-vanishing of its curvatures, one needs to assume that no
one of the four constants $a_1$, $a_2$, $\alpha_1$, $\alpha_2$ is
vanishing and, moreover, $\alpha_1^2\not= \alpha_2^2$. If some of
these constants are negative, then one may apply a symmetry
transformation in $\mathbb R^4$ to make them positive. Thus, from
now on we will assume that $a_1$, $a_2$, $\alpha_1$, $\alpha_2$
are positive and, besides, $\alpha_1$ does not coincide with
$\alpha_2$.

Evidently, the curve $\gamma$ represented by
(\ref{position_vector}) is situated in a Clifford torus
$T^2\subset \mathbb R^4$ given implicitly by $x_1^2+x_2^2=a_1^2$,
$x_3^2+x_4^2=a_2^2$. The curve is closed if and only if
$\frac{\alpha_1}{\alpha_2}\in\mathbb Q$, i.e.,
$\frac{\alpha_1}{\alpha_2}=\frac{m_1}{m_2}$, where $m_1$ and $m_2$
are coprime integers. The minimal period $T$ for the parameter $t$
is expressed by the obvious formulae
\begin{eqnarray}\label{T}
T=\frac{2\pi m_1}{\alpha_1}=\frac{2\pi m_2}{\alpha_2}.
\end{eqnarray}

Our aim is to analyze the value of $\int\limits_\gamma
\sqrt{k_{1}^2+k_{2}^2+k_{3}^2} ds $. Substituting
(\ref{length})-(\ref{k3}) and taking into account (\ref{T}), one
gets the following.

\medskip

{\bf Proposition 4.}
\begin{eqnarray}\label{Pr4}
\int\limits_\gamma \sqrt{k_{1}^2+k_{2}^2+k_{3}^2} ds \,\, =\,\,
 2\pi\, \sqrt{m_1^2+m_2^2}.
\end{eqnarray}

We would underline that an arc-length and the curvatures of
$\gamma$ depend on all the constants $a_1$, $a_2$, $\alpha_1$,
$\alpha_2$. However, the value of $\int\limits_\gamma
\sqrt{k_{1}^2+k_{2}^2+k_{3}^2} ds$ depends on the coprime integers
$m_1$ and $m_2$ only! Consequently, the range of possible values
of the integral in question is countable.

The minimal possible value of the integral is equal to
$2\sqrt{5}\pi$, it is achieved if either $m_1=1$, $m_2=2$ or
$m_1=2$, $m_2=1$.\footnote{The cases of $(m_1,m_2)$ equal to
$(0,0)$, $(1,0)$, $(0,1)$, $(1,1)$,  which give to $
\sqrt{m_1^2+m_2^2}$ values less than $\sqrt{5}$, are prohibited
because $\alpha_1$ and $\alpha_2$ are positive and don't coincide.
} This completes the proof of Theorem 1.

Remark. The coprime integers $m_1$ and $m_2$ are of a topological
nature, since they are clearly relied to the fundamental group of
the torus $T^2$. Therefore, one may conjecture that the values of
$\int\limits_\gamma \sqrt{k_{1}^2+k_{2}^2+k_{3}^2} ds$ in the
general case of curves with non-constant curvatures have to obey
some kinds of topological (homotopical) restrictions.

\bigskip

\bigskip

{\bf \large 2. Osculating indicatrices of closed curves of
constant curvatures in $\mathbb R^4$}

\medskip

Now let construct the osculating indicatrix of the curve $\gamma$,
c.f. [4]. By definition, the osculating plane of $\gamma$ at an
arbitrary point $p\in \gamma$ is spanned and oriented by the first
two vectors of the Frenet fame of $\gamma$ at $p$. Equivalently,
it is spanned and oriented by the vectors $\frac{dx}{dt}(t)$,
$\frac{d^2x}{dt^2}(t)$, where $x=x(t)$ is the position-vector of
$p\in\gamma$. Being translated to the origin $O\in \mathbb R^4$,
this osculating plane represents a point in the Grassmann
manifolds $G(2,4)$\footnote{For definitions and geometrical
properties of $G(2,4)$ see [5, Ch.8], [6].}. By moving $p$ along
$\gamma$, one obtains the one-dimensional family of osculating
planes of $\gamma$, which generates a curve $\tilde \gamma$ in
$G(2,4)$. This curve is called the {\it osculating indicatrix} of
$\gamma$.

The Grassmann manifold $G(2,4)$ is embedded into $\mathbb R^6$
with the help of Plucker coordinates, see [5,Ch.8.2], [6].
Consequently, $\tilde \gamma \subset G(2,4)$ may be viewed as a
curve in $\mathbb R^6$. If $x=x(t)$ is the position-vector of
$\gamma$, then $\tilde \gamma$ is represented in $\mathbb R^6$ by
the position-vector
\begin{equation}\label{oi}
\tilde x = \frac{ \left[ \frac{dx}{dt}, \frac{d^2x}{dt^2}\right]}
{\vert\left[ \frac{dx}{dt}, \frac{d^2x}{dt^2}\right]\vert},
\end{equation}
where the brackets $\left[,\right]$ denote the exterior product of
vectors. More precisely, one has
\begin{equation}\label{oi2}
\left[ \frac{dx}{dt}, \frac{d^2x}{dt^2}\right] = \left(\tilde
x_{12}, \tilde x_{13}, \tilde x_{14}, \tilde x_{23}, \tilde
x_{24}, \tilde x_{34}\right), \end{equation} where
\begin{equation}\label{oi3}\tilde x_{ij} =
\frac{dx_i}{dt} \frac{d^2x_j}{dt^2} -
\frac{dx_j}{dt}\frac{d^2x_i}{dt^2}, \quad 1\leq i<j\leq 4.
\end{equation}

Recall that the position-vector $x(t)$ of $\gamma$ is given by
(\ref{position_vector}). Substituting into (\ref{oi})-(\ref{oi3}),
we can easily derive the position-vector $\tilde x(t)$ of the
osculating indicatrix $\tilde\gamma$ in $\mathbb R^6$:
\begin{equation}\label{oipv}
\tilde x = \left( \frac{a_1^2\alpha_1^3} {\lambda},
\frac{-a_1a_2\alpha_1\alpha_2} {\lambda}\left(
\alpha_1\cos\alpha_1t\sin\alpha_2t -
\alpha_2\sin\alpha_1t\cos\alpha_2t\right), \right.
\end{equation}
\begin{equation*}
\frac{-a_1a_2\alpha_1\alpha_2} {\lambda}\left(
\alpha_1\cos\alpha_1t\cos\alpha_2t +
\alpha_2\sin\alpha_1t\sin\alpha_2t\right),
\end{equation*}
\begin{equation*}
\frac{a_1a_2\alpha_1\alpha_2} {\lambda}\left(
\alpha_1\sin\alpha_1t\sin\alpha_2t +
\alpha_2\cos\alpha_1t\cos\alpha_2t\right),
\end{equation*}
\begin{equation*} \left.
\frac{a_1a_2\alpha_1\alpha_2} {\lambda}\left(
\alpha_1\sin\alpha_1t\cos\alpha_2t -
\alpha_2\cos\alpha_1t\sin\alpha_2t\right), \frac{-a_2^2\alpha_2^3}
{\lambda}
 \right) ,
\end{equation*}
where $\lambda=\sqrt{a_1^2\alpha_1^2+a_2^2\alpha_2^2}
\sqrt{a_1^2\alpha_1^4+a_2^2\alpha_2^4}$.

Clearly, the curve $\tilde\gamma$ belongs to the four-dimensional
affine plane $\mathbb R^4\subset \mathbb R^6$ given by the
equations $\tilde x^{12}=\frac{a_1^2\alpha_1^3} {\lambda}$,
$\tilde x^{34}=\frac{-a_2^2\alpha_2^3} {\lambda}$.

By using standard formulae of the classical theory of curves, see
[1,Ch.32], one may find from (\ref{oipv}) the following
expressions for an arc length $\tilde s$ and the curvatures
$\tilde k_1$, $\tilde k_2$, $\tilde k_3$ of $\tilde \gamma$.

\medskip

{\bf Proposition 5.} {\it

1. An arc length $\tilde s$  of $\tilde\gamma$ is expressed as
follows:
\begin{eqnarray}\label{ts}
\frac{d\tilde s}{dt} =\frac{1}{\lambda}
\frac{a_1a_2\alpha_1\alpha_2}{ \vert\alpha_1^2-\alpha_2^2\vert}.
\end{eqnarray}

2. The curvatures $\tilde k_1$, $\tilde k_2$, $\tilde k_3$ of
$\tilde\gamma$ are expressed as follows:
\begin{eqnarray}\label{tk1}
\tilde k_1 = \lambda
\frac{\sqrt{\alpha_1^2+\alpha_2^2}}{a_1a_2\alpha_1\alpha_2\vert\alpha_1^2-\alpha_2^2\vert},
\end{eqnarray}
\begin{eqnarray}\label{tk2}
\tilde k_2 =\lambda
\frac{2}{a_1a_2\vert\alpha_1^2-\alpha_2^2\vert\sqrt{\alpha_1^2+\alpha_2^2}},
\end{eqnarray}
\begin{eqnarray}\label{tk3} \tilde k_3 = \lambda
\frac{1}{a_1a_2\alpha_1\alpha_2\sqrt{\alpha_1^2+\alpha_2^2}}.
\end{eqnarray}}

Consequently, the curvatures $\tilde k_1$, $\tilde k_2$, $\tilde
k_3$ of $\tilde \gamma$ are constant. Moreover, by applying
(\ref{ts})-(\ref{tk3}), one can easily verify that
(\ref{4_4})-(\ref{equality}) hold true, and this completes the
proof of Theorem 2.

\begin{small} Notice that if $\alpha_1=\alpha_2$ then $\gamma$ represented by
(\ref{position_vector}) is a circle, its first curvature is
$k_1=\frac{1}{\sqrt{a_1^2+a_2^2}}$, the second curvature is
vanishing, $k_2=0$, and the third curvature $k_3$ is undefined. In
this case all the osculating planes of $\gamma$ coincide with the
two-dimensional plane containing $\gamma$. Hence the osculating
indicatrix $\tilde\gamma$ degenerates to a point in
$G(2,4)$.\end{small}

\bigskip

\bigskip

{\bf \large 3. Concluding remarks and questions}

\medskip

{\it Remark  1.} Theorem 1 may be extended to the case of closed
curves of constant curvatures in $E^{2n}$, $n>2$. Moreover, one
can consider the integrals $\int\limits_\gamma
\sqrt{k_{2m-1}^2+k_{2m}^2+k_{2m+1}^2} ds$ for closed curves with
non-vanishing constant curvatures in $E^{2n}$, $1\leq m\leq n-1$.
It turns out that these integrals satisfy the same estimate
(\ref{4_3}), although the proof is more complicated technically.

{\it Remark  2.} Let $\gamma$ be an arbitrary smooth curve with
non-vanishing curvatures $k_1$, $k_2$ in $\mathbb R^4$. Then its
osculating indicatrix $\tilde \gamma$ is a smooth curve in the
Grassmann manifold $G(2,4)$. By considering the Plucker
coordinates, one may embed $G(2,4)$ into the unit sphere
$S^5\subset \mathbb R^6$, see [5, Ch.8.2], [6]. It turns out that
the osculating indicatrix $\tilde \gamma$ is an {\it asymptotic}
curve of $G(2,4)\subset S^5$. The integral $\int\limits_\gamma
\sqrt{k_{1}^2+k_{2}^2+k_{3}^2} ds$ is equal to the total curvature
$\int\limits_{\tilde\gamma}\tilde k_1 d\tilde s$ of $\tilde
\gamma$, when $\tilde \gamma$ is viewed as a curve in $\mathbb
R^6$, hence the optimality problem for $\int\limits_\gamma
\sqrt{k_{1}^2+k_{2}^2+k_{3}^2} ds$ gives rise to the following
problem. Considering arbitrary closed asymptotic curves
$\tilde\gamma$ in $G(2,4)\subset S^5$, find an optimal estimate
for $\int\limits_{\tilde\gamma}\tilde k_1 d\tilde s$.

{\it Remark  3.} From a local point of view, a smooth curve
$\tilde \gamma\subset G(2,4)$ is locally the osculating indicatrix
of some smooth curve $\gamma\subset \mathbb R^4$ if and only if
$\tilde\gamma$ is an asymptotic curve in $G(2,4)\subset S^5$, c.f.
[7]. We are interested in a global version of this statement.
Namely, what are the necessary and sufficient conditions for a
{\it closed} smooth curve in $G(2,4)$ to be the osculating
indicatrix of a smooth {\it closed} curve in $\mathbb R^6$?

\section*{References}

\noindent [1] Aminov Yu.A., Differential geometry and topology of
curves -- CRC Press, 2001, 216 pp.

\noindent [2] Fenchel W., Uber Krummung und Windung geschlossener
Raumkurven
// Math. Ann. 101 (1929), 238-252.

\noindent [3] Borsuk K., Sur la courbure totale des courbes
ferm\'ees
// Annales de la Soc. Polonaise 20 (1947), 251-265.

\noindent [4] Gorkaviy V., One integral inequality for closed
curves in Euclidean space // C. R. Acad. Sci. Paris 321 (1995),
1587-1591.

\noindent [5] Aminov Yu.A., Geometry of submanifolds -- CRC Press,
2001, 371 pp.

\noindent [6] Borisenko A.A., Nikolaevskii Yu.A., Grassmann
manifolds and the Grassmann image of submanifolds // Russian
Mathematical Surveys 46:2 (1991), 45–94.

\noindent [7] Gor'kavyi V., Reconstruction of a submanifold of
Euclidean space from its Grassmannian image that degenerates into
a line // Mathematical Notes 59:5 (1996), 490–497.

\bigskip

\bigskip

Vasyl Gorkavyy

B.Verkin Institute for Low Temperature Physics and Engineering,

47 Nauka ave., 61103 Kharkiv, Ukraine

{\it e-mail:} gorkaviy@ilt.kharkov.ua

\medskip

Raisa Posylaieva

Kharkiv National University of Civil Engineering and Architecture,

40 Sumska str.,  61002 Kharkiv, Ukraine

\end{document}